\documentclass[graybox]{svmult}

\usepackage{type1cm}        
\usepackage{makeidx}         
\usepackage{graphicx}        
\usepackage{multicol}        
\usepackage[bottom]{footmisc}

\usepackage{newtxtext}       %
\usepackage[varvw]{newtxmath}       

\usepackage{arydshln}

\usepackage{tikz,pgfplots}
\usetikzlibrary{intersections, shapes}
\usetikzlibrary{arrows}

\makeindex             

\begin{document}

\title*{Dual-Primal Isogeometric Tearing and Interconnecting methods for the Stokes problem}
\titlerunning{IETI-DP for the Stokes problem}
\author{Jarle Sogn \inst{*} and Stefan Takacs}
\institute{Jarle Sogn \at Johann Radon Institute for Computational and Applied Mathematics,
Austrian Academy of Sciences, Altenberger Str. 69, 4040 Linz, Austria \email{jarle.sogn@ricam.oeaw.ac.at}
\and Stefan Takacs \at Institute of Computational Mathematics,
Johannes Kepler University Linz, Altenberger Str. 69, 4040 Linz, Austria \email{stefan.takacs@numa.uni-linz.ac.at}
\\[.5em] \inst{*} Corresponding Author}
%
%
\maketitle

\abstract{
  We are interested in a fast solver for linear systems obtained by discretizing the Stokes problem with multi-patch Isogeometric Analysis.
  We use Dual-Primal Isogeometric Tearing and Interconnecting (IETI-DP) methods.
  In resent years, IETI-DP and related methods have been studied extensively, mainly for the Poisson problem. 
  For the Stokes equations, several challenges arise since the corresponding system is not positive definite, but has saddle point structure.
  Moreover, the Stokes equations with Dirichlet boundary conditions have a null-space, consisting of the constant pressure modes. This poses a challenge when considering the scaled Dirichlet preconditioner.
  We test out two different scaled Dirichlet preconditioners with different choices of primal degrees of freedom.
  The tests are performed on rather simple domains (the unit square and a quarter annulus) and a more complicated domain (a Yeti-footprint).}

\section{Introduction}
\label{SognTakacs:sec:intro}
We explore fast solvers for linear systems that arise from the discretization of the Stokes problem using Isogeometric Analysis (IgA; \cite{hughes2005isogeometric, da2014mathematical}).
We consider multipatch domains, that is, the computational domains consists of multiple non-overlapping patches. For such domains, FETI-DP methods (introduced in \cite{farhat2001feti}) are a canonical choice. FETI-DP was first adapted to IgA in \cite{kleiss2012ieti} and named the \textit{Dual-Primal Isogeometric Tearing and Interconnecting} (IETI-DP) method.
Several IETI-DP solvers for second-order elliptic boundary value problems have already been explored, see, e.g., \cite{hofer2017dual, hofer2019dual} and the more recent convergence analysis \cite{SchneckenleitnerTakacs:2020}, which is, besides grid sizes and the patch diameters, also robust in the spline degree and spline smoothness.

In~\cite{pavarino2016isogeometric}, a FETI-DP like solver has been applied to an isogeometric Taylor-Hood element for a single-patch domain. There, the substructures used for the setup of the solvers are non-overlapping parts of the considered patch. The isogeometric Taylor-Hood element uses spline degree $p+1$ for velocity and spline degree $p$ for pressure and smoothness $p-1$ for both of them. The derived solver maintains this smoothness. In~\cite{pavarino2016isogeometric}, a solver has also been proposed for the elasticity problem for quasi incompressible materials. These results have recently been extended in~\cite{widlund2021block}. In the finite element literature, there are a few results on FETI-DP methods for the Stokes problem, see, e.g., \cite{kimleepark,li2005dual,tuli} for the case of two dimensions and \cite{tu2015feti} for the case of three dimensions.

We consider IETI-DP solvers for multi-patch domains, where each patch is a substructure. Concerning the coupling between the patches, 
we go a slightly different way and use the minimum smoothness requirements that are feasible in order to obtain a conforming discretization. This means that the pressure space is discontinuous across the interfaces between the patches. For the velocity, we only impose continuity between the patches. Within the patches, we use the aforementioned isogeometric Taylor-Hood scheme.
Several challenges arise when one aims to extend the IETI-DP solvers to the Stokes problem. The first is that the resulting linear system is not positive definite, it rather has a saddle point structure. The second challenge is that in case of Dirichlet conditions, the differential operator has a nonzero nullspace. The nullspace consists of all constant pressure modes. This is normally remedied by restricting the solution space such that the average pressure is zero. When one wants to apply a scaled Dirichlet preconditioner, such a condition would be required locally. We obtain such a local condition by adding the patchwise averages of the pressure to the space of primal degrees of freedom.

The remainder of this paper is organized as follows. In Section~\ref{SognTakacs:sec:IgAandIETI}, we formulate the problem and introduce a few possible ways of realizing the IETI-DP solver. In Section~\ref{SognTakacs:sec:num}, we present a numerical results. Finally, we draw some conclusions in Section~\ref{SognTakacs:sec:conclusions}.

\section{Problem formulation and IETI-DP}
\label{SognTakacs:sec:IgAandIETI}
As model problem, we consider the Stokes problem: Let $\Omega\subset \mathbb{R}^2$ be a bounded Lipschitz domain. For a given right-hand side $\mathbf{f} \in \left[L(\Omega)\right]^2$, find $(\mathbf{u},p)\in \left[H^1_0(\Omega)\right]^2\times L^2_0(\Omega)$ such that
\begin{equation}
  \label{SognTakacs:eq:stokes}
  \begin{aligned} 
	(\nabla \mathbf{u}, \nabla \mathbf{v})_{L^2(\Omega)} + (p, \mathrm{div}\, \mathbf{v})_{L^2(\Omega)} &= (\mathbf{f},\mathbf{v})_{L^2(\Omega)} \quad &&\forall\, \mathbf{v} \in \left[H^1_0(\Omega)\right]^2,\\
	(\mathrm{div}\, \mathbf{u}, q)_{L^2(\Omega)} &= 0  \quad &&\forall\, q \in L^2_0(\Omega).
  \end{aligned} 
\end{equation}
We assume that the domain $\Omega$ is composed of $K$ non-overlapping patches $\Omega^{(k)}$ which are parameterized with geometry mappings $\mathbf{G}_k:(0,1)^2\rightarrow  \Omega^{(k)}$ such that
$\Omega^{(k)}= \mathbf{G}_k((0,1)^2)$. We assume that both the Jacobian $\nabla \mathbf{G}_k$ and its inverse $(\nabla \mathbf{G}_k)^{-1}$ are almost everywhere uniformly bounded (cf.~\cite[Ass.~1]{SchneckenleitnerTakacs:2020}).

As local discretization spaces on the parameter domain $(0,1)^2$, we use tensor-product B-splines which we denote by $S_{p,\alpha}$, where $p$ is the spline degree and $\alpha$ is the smoothness, this is, such that the splines are $\alpha$ times continuously differentiable. So, $\alpha = p-1$ corresponds to splines of maximum smoothness. For the Stokes equations, we need a inf-sup stable discretization space. We use the generalized Taylor-Hood space (cf. \cite{bressan2013isogeometric}), which we transfer to the physical domain $\Omega^{(k)}$ using the pull-back principle. The spaces on the individual patches are given by
\begin{equation}
\label{SognTakacs:eq:lccalspaces}
\begin{aligned}
  \mathbf{V}^{(k)} &:= \left\lbrace \mathbf{u}:
  		\mathbf{u}\circ \mathbf{G}_k \in [S_{p+1,\alpha}]^2
  		\mbox{ and }
  		\mathbf{u}|_{\partial\Omega\cap\partial\Omega^{(k)}} = 0
  		\right\rbrace,\\
  		Q^{(k)} &:= \{q : q\circ \mathbf{G}_k\in S_{p,\alpha} \},
\end{aligned}
\end{equation}
where $\mathbf{u}|_{\partial\Omega\cap\partial\Omega^{(k)}}$ denotes the restriction of $\mathbf{u}$ to $\partial\Omega\cap\partial\Omega^{(k)}$ (trace operator).

The overall discretization space for the velocity is $[H^1_0(\Omega)]^2$, therefore a conforming discretization needs to be continuous. Thus, we assume that the spaces are fully matching, that is, for each basis function active on an interface between to patches, there is exactly one basis function on the neighboring patch such that the functions agree on the interface (cf.~\cite[Ass.~5]{SchneckenleitnerTakacs:2020}). This allows us to define the overall discretization space for the velocity via $\mathbf{V}:= \{\mathbf u\in [H^1_0(\Omega)]^2 : \mathbf u|_{\Omega^{(k)}} \in \mathbf{V}^{(k)} \}$.
For the pressure, continuity is not required in order to obtain a conforming discretization. We define $Q:= \{q\in L^2_0(\Omega) : q|_{\Omega^{(k)}} \in Q \}$, where $L^2_0(\Omega)$ is the space of functions on $L^2(\Omega)$ with vanishing mean value.

By using the spaces from~\eqref{SognTakacs:eq:lccalspaces} for a patch-wise assembling of the problem~\eqref{SognTakacs:eq:stokes}, we obtain local (still uncoupled) linear
systems
\[
    A^{(k)}\underline{\mathbf{x}}^{(k)}=
		\begin{pmatrix}
		  K^{(k)} & (D^{(k)})^\top \\
		  D^{(k)} & 0 
		\end{pmatrix}
        \begin{pmatrix}
          \underline{\mathbf{u}}^{(k)}\\
          \underline{p}^{(k)}
        \end{pmatrix}
        =
        \begin{pmatrix}
          \underline{\mathbf{f}}^{(k)} \\
          0
        \end{pmatrix}
        = \underline{\mathbf{b}}^{(k)},
\]
where the matrices $K^{(k)}$ and $D^{(k)}$ and the vector $\underline{\mathbf{f}}^{(k)}$ are obtained by evaluating the corresponding terms $(\nabla \cdot, \nabla \cdot)_{L^2(\Omega^{(k)})}$, $(\cdot, \mathrm{div}\, \cdot)_{L^2(\Omega^{(k)})}$ and $(\mathbf{f},\cdot)_{L^2(\Omega^{(k)})}$, respectively, for the basis functions of the bases for the spaces $\mathbf{V}^{(k)}$ and/or $Q^{(k)}$.
An underlined quantity, like $\underline{\mathbf{u}}^{(k)}$, represents the coefficient vector of the corresponding function ${\mathbf{u}}^{(k)}$ with respect to the chosen basis.

Next, we split the degrees of freedom for the velocity variable into those associated to the interfaces, denoted $\underline{\mathbf{u}}_\Gamma$, and the remaining degrees of freedom, associated to basis functions vanishing on the interface, denoted $\underline{\mathbf{u}}_{\text{I}}$. Using this splitting, we obtain
\begin{equation}
  \nonumber
    A^{(k)}\underline{\mathbf{x}}^{(k)}=
		\begin{pmatrix}
		  K_{\Gamma\Gamma}^{(k)} &K_{\Gamma\text{I}}^{(k)} & (D_{\Gamma}^{(k)})^\top \\
		  K_{\text{I}\Gamma}^{(k)}     &K_{\text{II}}^{(k)}     & (D_{\text{I}}^{(k)})^\top \\
		  D_\Gamma^{(k)}         &D_{\text{I}}^{(k)}         & 0 
		\end{pmatrix}
        \begin{pmatrix}
          \underline{\mathbf{u}}_\Gamma^{(k)}\\
          \underline{\mathbf{u}}_{\text{I}}^{(k)}\\
          \underline{p}^{(k)}
        \end{pmatrix}
        =
        \begin{pmatrix}
          \underline{\mathbf{f}}_\Gamma^{(k)} \\
          \underline{\mathbf{f}}_{\text{I}}^{(k)} \\
          0
        \end{pmatrix}
        = \underline{\mathbf{b}}^{(k)}.
\end{equation}
Next, we introduce primal degrees of freedom, which are introduced in order to guarantee that the system matrices of the patch-local problems are non-singular. As in the case of the Poisson problem, the patch-local systems of patches that do not contribute to the Dirichlet boundary represent problems with pure Neumann boundary conditions. On these patches, the null space of $A^{(k)}$ corresponds to constant velocity modes. We mitigate this problem by strongly enforcing the continuity for the primal degrees of freedom. We consider three variants:
\begin{itemize}
\item Variant $\Pi_{\mathrm{c}}$: Each velocity component on each of the corners.
\item Variant $\Pi_{\mathrm{ce}}$: Each velocity component on each of the corners and the average of each velocity component on each of the edges.
\item Variant $\Pi_{\mathrm{cn}}$: Each velocity component on each of the corners and the average of the normal component of the velocity on each of the edges.
\end{itemize}
In each of the local problems, we require that these primal degrees of freedom vanish. This constraint is represented by the condition $C_v^{(k)} \underline{\mathbf u}_\Gamma^{(k)}=0$.

Moreover, we introduce primal degrees of freedom for the pressure. For the scaled Dirichlet preconditioner, it is necessary that the local Dirichlet problems are uniquely solvable. However, the Stokes system with Dirichlet boundary conditions has a non-trivial null space, which consists of the constant pressure modes. We also choose the averages of the pressure on the individual patches as primal degrees of freedom. For each of the local problems we require that this primal degree of freedom vanishes. This constraint is represented by the condition $C_p^{(k)} \underline{p}^{(k)}=0$.

The local systems read as follows
\begin{equation}
  \label{SognTakacs:eq:localStokesAverage}
  \begin{aligned}
    \bar{A}^{(k)}\bar{\underline{\mathbf{x}}}^{(k)}
    &=
    \begin{pmatrix}
    		A^{(k)} & (C^{(k)})^\top \\
    		C^{(k)} & 0
		\end{pmatrix}
		\begin{pmatrix}
				\underline{\mathbf{x}}^{(k)} \\
				\underline{\mu}^{(k)}
		\end{pmatrix}
		\\&=
		\left(
		\begin{array}{ccc;{.5pt/1pt}cc}
		  K_{\Gamma\Gamma}^{(k)} &K_{\Gamma\text{I}}^{(k)}  & (D_{\Gamma}^{(k)})^\top &0&C_v^\top\\
		  K_{\text{I}\Gamma}^{(k)}     &K_{\text{I}\text{I}}^{(k)}      & (D_{\text{I}}^{(k)})^\top     &0&0\\
		  D_\Gamma^{(k)}         &D_{\text{I}}^{(k)}          & 0               &(C_p^{(k)})^\top&0\\ \hdashline[.5pt/1pt]
		  0                &0             & C_p^{(k)}             &0&0\\
		  C_v^{(k)}              &0             & 0               &0&0\\		
		  \end{array}
		  \right)
        \begin{pmatrix}
          \underline{\mathbf{u}}_\Gamma^{(k)}\\
          \underline{\mathbf{u}}_{\text{I}}^{(k)}\\
          \underline{p}^{(k)}\\ \hdashline[.5pt/1pt]
          \mu_p^{(k)}\\
          \underline{\mu}_v^{(k)}
        \end{pmatrix}
        =
        \begin{pmatrix}
          \underline{\mathbf{f}}_\Gamma^{(k)} \\
          \underline{\mathbf{f}}_{\text{I}}^{(k)} \\
          0\\\hdashline[.5pt/1pt]
          0\\
          0
        \end{pmatrix}=\bar{\underline{\mathbf{b}}}^{(k)},
   \end{aligned}
\end{equation}
where $\mu_p^{(k)}$ and $\underline{\mu}_v^{(k)}$ are the Lagrangian multipliers corresponding to the primal degrees of freedom.

The continuity of the velocity variables (except the corner values) between the patches is enforced weakly using a constraint that is represented by the boolean matrix $\bar{B}^{(k)}=(B_\Gamma^{(k)},0,0,0,0)$. The velocity is continuous if $\sum_{k=1}^K\bar{B}^{(k)} \bar{\underline{\mathbf{x}}}^{(k)} = 0$ or, equivalently, $\sum_{k=1}^K B_\Gamma^{(k)} \underline{\mathbf{u}}_\Gamma^{(k)} = 0$.

Finally, we consider the primal problem, this is, the global problem for the primal degrees of freedom. We use a nodal basis for the primal degrees of freedom which is $A^{(k)}$-orthogonal to the remaining degrees of freedom on each patch. This basis is represented with respect to the patch-local bases by the matrices $\Psi^{(k)}$, which are the solutions to
\[
    \begin{pmatrix}
    		A^{(k)} & (C^{(k)})^\top \\
    		C^{(k)} & 0
		\end{pmatrix}
		\begin{pmatrix}
			\Psi^{(k)} \\ \mathrm{M}^{(k)}
		\end{pmatrix}
		=
		\begin{pmatrix}
				0 \\ R_C^{(k)}
		\end{pmatrix},
\]
where $R_C^{(k)}$ is a boolean matrix, which represents the assignment between a patch-local ordering of the primal degrees of freedom and the global ordering of the primal degrees of freedom. Using the matrices $\Psi^{(k)}$, we define the system matrix, jump matrix and right-hand side for the primal problem as
\begin{equation*}
A_\Pi := \sum^K_{k=1} (\Psi^{(k)})^\top {A}^{(k)}\Psi^{(k)}, \quad B_\Pi := \sum^K_{k=1} {B}^{(k)}\Psi^{(k)} \quad \text{and}\quad \mathbf{b}_\Pi := \sum^K_{k=1} (\Psi^{(k)})^\top \underline{{\mathbf{b}}}^{(k)}.
\end{equation*}
Before we can write up the global system, we need to revisit the condition that the average of the pressure vanishes, this is, $p\in L^2_0(\Omega)$. Our choice of the primal degrees of freedom ensures that the average of the pressure vanishes for all patch-local problems. Since the patch-wise constant pressure modes are primal, they form part of the primal system. In order to obtain the unique solvability of the primal system, we introduce a constant that guarantees that the global average of the pressure vanishes. This condition is enforced by adding a vector $\mathbf{p}$ representing this constraint. So, we extend the primal system as follows 
\begin{equation}\label{SognTakacs:eq:Apibar}
\bar{A}_\Pi := 
\begin{pmatrix}
  A_\Pi &\mathbf{p}^\top\\
  \mathbf{p} & 0
\end{pmatrix},
\quad \bar{B}_\Pi := \big(B_\Pi,\, 0\big),
  \quad
  \bar{\underline{\mathbf{x}}}_\Pi :=
  \begin{pmatrix}
  \underline{\mathbf{x}}_\Pi \\
  \mu_0
  \end{pmatrix}\quad\text{and} \quad \bar{\underline{\mathbf{b}}}_\Pi :=
\begin{pmatrix}
  \underline{\mathbf{b}}_\Pi \\
  0 \end{pmatrix},
\end{equation}
where $\mu_0$ is a new Lagrangian multiplier associated to the additional constraint.

By coupling the patch local systems and the primal system, we obtain the IETI-DP saddle point system, which characterizes the solution and which reads as follows
\begin{equation}\label{SognTakacs:eq:bigspp}
\begin{pmatrix}
	\bar A^{(1)}&        &              &              & (\bar B^{(1)})^\top \\
	            & \ddots &              &              & \vdots  \\
	            &        & \bar A^{(K)} &              & (\bar B^{(K)})^\top \\
              &        &              & \bar A_\Pi   & \bar B_\Pi^\top \\
  \bar B^{(1)}& \cdots & \bar B^{(K)} & \bar B_\Pi   & 0
\end{pmatrix}
\begin{pmatrix}
	\underline{\bar{\mathbf{x}}}^{(1)} \\
	\vdots \\
	\underline{\bar{\mathbf{x}}}^{(K)} \\
	\underline{\bar{\mathbf{x}}}_\Pi \\
	\underline{\lambda}
\end{pmatrix}
=
\begin{pmatrix}
	\underline{\bar{\mathbf{b}}}^{(1)} \\
	\vdots \\
	\underline{\bar{\mathbf{b}}}^{(K)} \\
	\underline{\bar{\mathbf{b}}}_\Pi \\
	0
\end{pmatrix}
\end{equation}
We solve this linear system as follows. We first define the Schur-complement $\bar{F}$ and the corresponding right-hand side $\underline g$ by
\begin{equation*}
  \bar{F} := \bar{B}_\Pi \bar{A}^{-1}_\Pi \bar{B}_\Pi^\top+\sum^{K}_{k = 1} \bar{B}^{(k)} (\bar{A}^{(k)})^{-1}(\bar{B}^{(k)})^\top,\; \underline{g}:= \bar{B}_\Pi \bar{A}^{-1}_\Pi\bar{\underline{\mathbf{b}}}_\Pi + \sum^{K}_{k = 1}  \bar{B}^{(k)} (\bar{A}^{(k)})^{-1}\bar{\underline{\mathbf{b}}}^{(k)},
\end{equation*}
where we note that the symmetric and positive definite matrix $\bar F$ is never computed. Products of the form $\bar F \underline \lambda$ can be easily computed by solving the primal problem and the patch-local problems with right-hand sides
$\bar B_\Pi^\top \underline \lambda$ and $(\bar B^{(k)})^\top \underline \lambda$, respectively. Consequently, we solve the linear system
\begin{equation}
  \label{SognTakacs:eq:lambdaSys}
  \bar{F}\underline{\lambda} = \underline{g}
\end{equation}
using a preconditioned conjugate gradient solver. After the computation of $\underline{\lambda}$, the local solutions are recovered using~\eqref{SognTakacs:eq:bigspp}:
\[
		\bar{\underline{\mathbf{x}}}^{(k)} = (\bar A^{(k)})^{-1}
		(\bar{\underline{\mathbf{b}}}^{(k)}
		-(\bar{B}^{(k)})^\top \underline\lambda)
		\quad\mbox{and}\quad
		\bar{\underline{\mathbf{x}}}_\Pi = \bar A_\Pi^{-1}
		(\bar{\underline{\mathbf{b}}}_\Pi
		-\bar{B}_\Pi^\top \underline\lambda).
\]
The local solution $\underline{\mathbf{x}}^{(k)}=(\underline{\mathbf{u}}^{(k)},\underline{p}^{(k)})$ is then obtained by
$
			\underline{\mathbf{x}}^{(k)} =
			\big( I,\,0 \big)\; \bar{\underline{\mathbf{x}}}^{(k)}
			+
			\Psi^{(k)} \big( I,\,0 \big)\; \bar{\underline{\mathbf{x}}}_\Pi,
$
where the matrices $(I,0)$ eliminate the rows corresponding to the Lagrangian multipliers $\mu_p$, $\underline{\mu}_v$ and $\mu_0$, cf.~\eqref{SognTakacs:eq:localStokesAverage} and
\eqref{SognTakacs:eq:Apibar}.

Finally, we have to discuss the choice of the preconditioner for the conjugate gradient solver. We consider two variants of a scaled  Dirichlet preconditioner. The first variant follows the standard construction principles and is based on the equation of interest, which is the Stokes equation. We choose
\begin{equation}\nonumber
		M_{\mathrm{sD},1} := 
		\sum_{k=1}^K B_{\Gamma}^{(k)} (D^{(k)})^{-1} S_{\Gamma,1}^{(k)} (D^{(k)})^{-1} (B_{\Gamma}^{(k)})^\top,
\end{equation}
where the corresponding Schur complement is given by
\begin{equation*}
		S_{\Gamma, 1}^{(k)} :=
		K_{\Gamma\Gamma}^{(k)}-
		\begin{pmatrix}
				K_{\Gamma\text{I}}^{(k)}         & (D_\Gamma^{(k)})^\top        & 0
		\end{pmatrix}
		\begin{pmatrix}
				K_{\text{I}\text{I}}^{(k)}       & (D_{\text{I}}^{(k)})^\top & 0\\
				D_{\text{I}}^{(k)}           & 0           & (C_p^{(k)})^\top \\
				0              & C_p^{(k)}         & 0
		\end{pmatrix}^{-1}
		\begin{pmatrix}
				K_{\text{I}\Gamma}^{(k)} \\
				D_{\Gamma}^{(k)}   \\
				0
		\end{pmatrix}
\end{equation*}
and
$D^{(k)}:=2I$ is set up based on the principle of multiplicity scaling. Often, the primal degrees of freedom do not affect the definition of the Schur complement used for the scaled Dirichlet preconditioner. Note that our definition includes the constraint that guarantees that the pressure average vanishes. This is required in order to make the local system matrix solvable.

The second variant is motivated by the observation that the scaled Dirichlet preconditioner is usually set up in order to realize a $H^{1/2}$-scalar product. To archive this goal, we simply choose
\begin{equation}\nonumber
		M_{\mathrm{sD},2} := 
		\sum_{k=1}^K B_{\Gamma}^{(k)} (D^{(k)})^{-1} S_{\Gamma,2}^{(k)} (D^{(k)})^{-1} (B_{\Gamma}^{(k)})^\top,
\end{equation}
where
\[
	S_{\Gamma, 2}^{(k)} := K_{\Gamma\Gamma}^{(k)} - K_{\Gamma\text{I}}^{(k)} (K_{\text{I}\text{I}}^{(k)})^{-1} K_{\text{I}\Gamma}^{(k)},
\]
which corresponds to a vector valued Poisson problem.

\section{Numerical experiments}
\label{SognTakacs:sec:num}
\begin{figure}[htb]
	\centering
	\includegraphics[height=4cm]{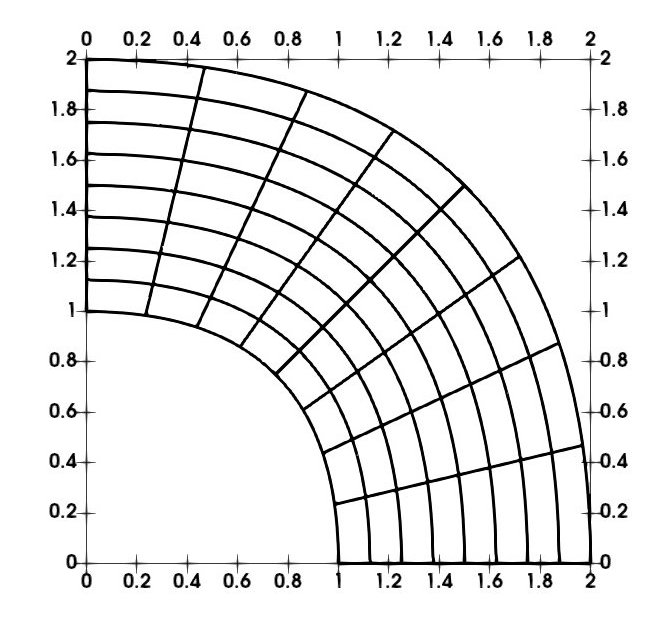}
    \qquad
	\includegraphics[height=4cm]{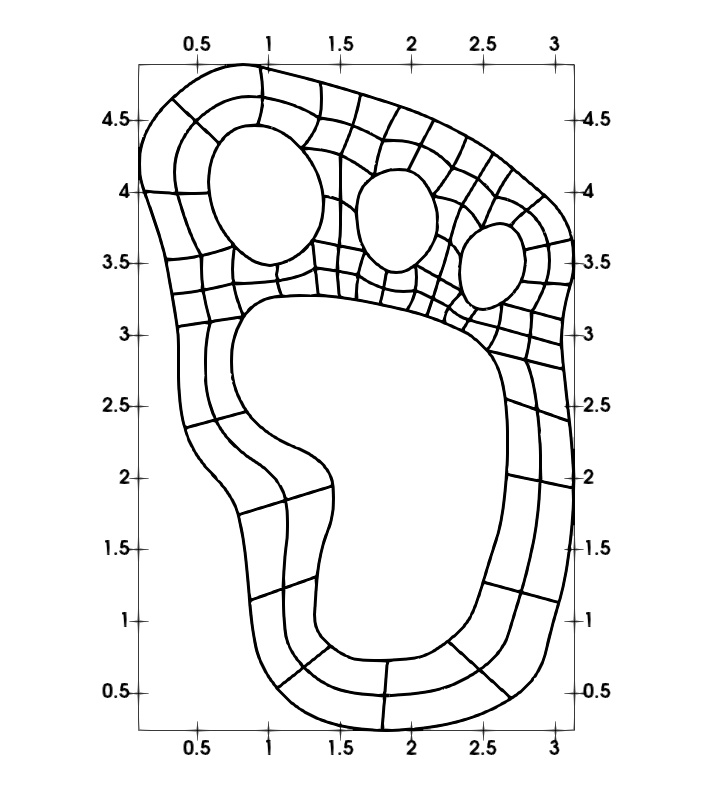}
    \qquad
	\raisebox{1.25cm}{
	\includegraphics[height=1.5cm]{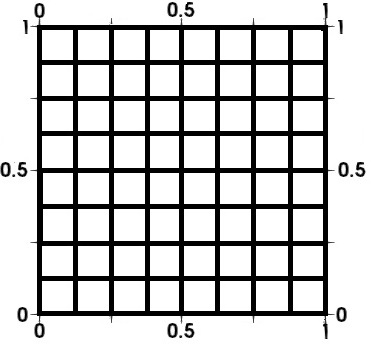}
	}
	\caption{Computational domains: Quarter annulus (left); Yeti-footprint (middle); unit square (right)}
	\label{SognTakacs:fig:domain}
\end{figure}

We consider the Stokes problem~\eqref{SognTakacs:eq:stokes} with the right-hand side function
\begin{align*}
  \mathbf{f}(x,y) &= (-\pi\cos(\pi x)-2 \pi^2\sin(\pi x) \cos(\pi y),\, 2\pi^2\cos(\pi x) \sin(\pi y))
\end{align*}
and the inhomogeneous Dirichlet boundary conditions
\begin{align*}
  \mathbf{u}(x,y) &= (-\sin(\pi x)\cos(\pi y),\, \cos(\pi x)\sin(\pi y)) \quad \text{for} \quad (x,y)\in \partial\Omega.
\end{align*}
We consider three computational domains: a quarter annulus (64 patches), the Yeti-footprint (84 patches) and the unit square (64 patches), see Figure~\ref{SognTakacs:fig:domain}.
For the quarter annulus and the unit square, the patch-local function spaces on the coarsest grid level ($\ell=0$) consist of global polynomials only, this is, there are no inner knots. For the Yeti-footprint, the patch-local function spaces on the long and thin patches at the bottom consist of two elements, obtained by uniformly refining the two longer sides of the patch. For the remaining patches of the Yeti-footprint, the function spaces consist of global polynomials only. The grid levels $\ell=1,2,\ldots$ are obtained by $\ell$ uniform refinement levels. In either case, the grid size on the parameter domain is $\widehat{h}= 2^{-\ell}$. Within each patch, we have splines of degree $p+1$ and smoothness $p-1$ for the velocity and splines of degree $p$ and smoothness $p-1$ for the pressure. As outlined in the previous section, the coupling between the patches is continuous for the velocity and discontinuous for the pressure. 

In the following, we discuss the convergence of the preconditioned conjugate gradient solver, used to solve~\eqref{SognTakacs:eq:lambdaSys}. We start the conjugate gradient solver with a random initial guess and stop the iteration when Euclidean norm of the residual vector is reduced by a factor of $10^{-6}$ compared to the Euclidean norm of the initial residual  vector. The local linear systems that need to be solved in the pre-processing steps (computation of $\underline g$) and during the main iteration are realized using a sparse LU-solver. All experiments have been implemented using the G+smo library\footnote{\url{https://github.com/gismo/gismo}} and have been performed on the Radon1 cluster\footnote{\url{https://www.ricam.oeaw.ac.at/hpc/}} in Linz.

\begin{table}[btp]
  \begin{center} 
  \begin{tabular}{|c|c|c|c|c|c|}
    \hline
    $\;\ell\setminus p\;$ & {\quad$2$\quad} & {\quad$3$\quad} & {\quad$4$\quad} & {\quad$5$\quad} & {\quad$6$\quad} \\
    \hline
    \hline
  $2$  & $43$ & $42$ & $43$ & $43$ & $38$  \\ \hline
	$3$  & $47$ & $48$ & $48$ & $48$ & $44$  \\ \hline
	$4$  & $53$ & $51$ & $55$ & $53$ & $52$  \\ \hline
	$5$  & $56$ & $60$ & $61$ & $58$ & $58$  \\
    \hline
  \end{tabular}
  $\quad$
  \begin{tabular}{|c|c|c|c|c|c|}
    \hline
    $\;\ell\setminus p\;$ & {\quad$2$\quad} & {\quad$3$\quad} & {\quad$4$\quad} & {\quad$5$\quad} & {\quad$6$\quad} \\
    \hline
    \hline
  $2$  & $28$ & $27$ & $28$ & $27$ & $25$  \\ \hline
	$3$  & $29$ & $31$ & $30$ & $30$ & $29$  \\ \hline
	$4$  & $34$ & $34$ & $35$ & $33$ & $32$  \\ \hline
	$5$  & $35$ & $37$ & $38$ & $36$ & $35$   \\
    \hline
  \end{tabular}
    \caption{Iteration numbers for quarter annulus using $\Pi_{\mathrm{c}}$ and $M_{\mathrm{sD},1}$ (left), $M_{\mathrm{sD},2}$ (right)}
    \label{SognTakacs:tab:APC}
  \end{center} 
\end{table}

\begin{table}[btp]
  \begin{center} 
  \begin{tabular}{|c|c|c|c|c|c|}
    \hline
    $\;\ell\setminus p\;$ & {\quad$2$\quad} & {\quad$3$\quad} & {\quad$4$\quad} & {\quad$5$\quad} & {\quad$6$\quad} \\
    \hline
    \hline
  $2$  & $16$ & $16$ & $16$ & $15$ & $15$  \\ \hline
	$3$  & $17$ & $17$ & $17$ & $17$ & $16$  \\ \hline
	$4$  & $18$ & $18$ & $18$ & $18$ & $18$  \\ \hline
	$5$  & $20$ & $20$ & $20$ & $19$ & $19$  \\
    \hline
  \end{tabular}
  $\quad$
  \begin{tabular}{|c|c|c|c|c|c|}
    \hline
    $\;\ell\setminus p\;$ & {\quad$2$\quad} & {\quad$3$\quad} & {\quad$4$\quad} & {\quad$5$\quad} & {\quad$6$\quad} \\
    \hline
    \hline
    $2$  & $\;11\;$ & $\;12\;$ & $\;11\;$ & $\;11\;$ & $\;11\;$  \\ \hline
	$3$  & $13$ & $13$ & $13$ & $12$ & $12$  \\ \hline
	$4$  & $14$ & $14$ & $14$ & $13$ & $13$  \\ \hline
	$5$  & $15$ & $15$ & $15$ & $15$ & $14$   \\
    \hline
  \end{tabular}

    \caption{Iteration numbers for quarter annulus using $\Pi_{\mathrm{ce}}$ and $M_{\mathrm{sD},1}$ (left), $M_{\mathrm{sD},2}$ (right)}
    \label{SognTakacs:tab:APE}
  \end{center} 
\end{table}
\begin{table}[tp]
  \begin{center} 
  \begin{tabular}{|c|c|c|c|c|c|}
    \hline
    $\;\ell\setminus p\;$ & {\quad$2$\quad} & {\quad$3$\quad} & {\quad$4$\quad} & {\quad$5$\quad} & {\quad$6$\quad} \\
    \hline
    \hline
  $2$  & $22$ & $23$ & $22$ & $22$ & $22$  \\ \hline
	$3$  & $25$ & $25$ & $25$ & $24$ & $24$  \\ \hline
	$4$  & $26$ & $27$ & $27$ & $26$ & $25$  \\ \hline
	$5$  & $29$ & $29$ & $28$ & $28$ & $27$   \\
    \hline
  \end{tabular}
  $\quad$
\begin{tabular}{|c|c|c|c|c|c|}
    \hline
    $\;\ell\setminus p\;$ & {\quad$2$\quad} & {\quad$3$\quad} & {\quad$4$\quad} & {\quad$5$\quad} & {\quad$6$\quad} \\
    \hline
    \hline
  $2$  & $17$ & $17$ & $17$ & $17$ & $16$  \\ \hline
	$3$  & $18$ & $19$ & $19$ & $18$ & $18$  \\ \hline
	$4$  & $20$ & $20$ & $20$ & $20$ & $19$  \\ \hline
	$5$  & $22$ & $22$ & $22$ & $21$ & $21$   \\
    \hline
  \end{tabular}
    \caption{Iteration numbers for quarter annulus using $\Pi_{\mathrm{cn}}$ and $M_{\mathrm{sD},1}$ (left), $M_{\mathrm{sD},2}$ (right)}
    \label{SognTakacs:tab:APN}
  \end{center} 
\end{table}
\begin{table}[tp]
  \begin{center} 
    \begin{tabular}{|c|c|c|c|c|c|}
    \hline
    $\;\ell\setminus p\;$ & {\quad$2$\quad} & {\quad$3$\quad} & {\quad$4$\quad} & {\quad$5$\quad} & {\quad$6$\quad} \\
    \hline
    \hline
  $2$  & $143$ & $150$ & $133$ & $132$ & $126$  \\ \hline
	$3$  & $172$ & $175$ & $146$ & $183$ & $176$  \\ \hline
	$4$  & $191$ & $201$ & $203$ & $151$ & $174$  \\ \hline
	$5$  & $223$ & $245$ & $238$ & $240$ & $234$  \\
    \hline
  \end{tabular}
   $\quad$
    \begin{tabular}{|c|c|c|c|c|c|}
    \hline
    $\;\ell\setminus p\;$ & {\quad$2$\quad} & {\quad$3$\quad} & {\quad$4$\quad} & {\quad$5$\quad} & {\quad$6$\quad} \\
    \hline
    \hline
  $2$  & $76$ & $78$ & $67$ & $70$ & $63$  \\ \hline
	$3$  & $85$ & $86$ & $73$ & $85$ & $78$  \\ \hline
	$4$  & $94$ & $97$ & $95$ & $74$ & $72$  \\ \hline
	$5$  & $101$ & $112$ & $106$ & $109$ & $103$   \\
    \hline
  \end{tabular}
    \caption{Iteration numbers for Yeti-footprint using $\Pi_{\mathrm{c}}$ and $M_{\mathrm{sD},1}$ (left), $M_{\mathrm{sD},2}$ (right)}
    \label{SognTakacs:tab:YPC}
  \end{center} 
\end{table}
\begin{table}[tp]
  \begin{center} 
      \begin{tabular}{|c|c|c|c|c|c|}
    \hline
    $\;\ell\setminus p\;$ & {\quad$2$\quad} & {\quad$3$\quad} & {\quad$4$\quad} & {\quad$5$\quad} & {\quad$6$\quad} \\
    \hline
    \hline
  $2$  & $21$ & $22$ & $22$ & $23$ & $23$  \\ \hline
	$3$  & $25$ & $27$ & $24$ & $28$ & $26$  \\ \hline
	$4$  & $29$ & $30$ & $30$ & $30$ & $31$  \\ \hline
	$5$  & $31$ & $34$ & $33$ & $33$ & $35$   \\
    \hline
  \end{tabular}
    $\quad$
    \begin{tabular}{|c|c|c|c|c|c|}
    \hline
    $\;\ell\setminus p\;$ & {\quad$2$\quad} & {\quad$3$\quad} & {\quad$4$\quad} & {\quad$5$\quad} & {\quad$6$\quad} \\
    \hline
    \hline
    $2$  & $\;14\;$ & $\;14\;$ & $\;14\;$ & $\;14\;$ & $\;13\;$  \\ \hline
	$3$  & $15$ & $16$ & $15$ & $16$ & $15$  \\ \hline
	$4$  & $17$ & $18$ & $16$ & $16$ & $17$  \\ \hline
	$5$  & $19$ & $19$ & $19$ & $18$ & $18$   \\
    \hline
  \end{tabular}
    \caption{Iteration numbers for Yeti-footprint using $\Pi_{\mathrm{ce}}$ and $M_{\mathrm{sD},1}$ (left), $M_{\mathrm{sD},2}$ (right)}
    \label{SognTakacs:tab:YPE}
  \end{center} 
\end{table}
\begin{table}[tp]
  \begin{center} 
   \begin{tabular}{|c|c|c|c|c|c|}
    \hline
    $\;\ell\setminus p\;$ & {\quad$2$\quad} & {\quad$3$\quad} & {\quad$4$\quad} & {\quad$5$\quad} & {\quad$6$\quad} \\
    \hline
    \hline
  $2$  & $20$ & $21$ & $22$ & $22$ & $23$  \\ \hline
	$3$  & $23$ & $25$ & $26$ & $26$ & $26$  \\ \hline
	$4$  & $27$ & $28$ & $28$ & $28$ & $29$  \\ \hline
	$5$  & $30$ & $31$ & $31$ & $31$ & $32$   \\
    \hline
  \end{tabular}
  $\quad$
    \begin{tabular}{|c|c|c|c|c|c|}
    \hline
    $\;\ell\setminus p\;$ & {\quad$2$\quad} & {\quad$3$\quad} & {\quad$4$\quad} & {\quad$5$\quad} & {\quad$6$\quad} \\
    \hline
    \hline
  $2$  & $16$ & $17$ & $16$ & $16$ & $16$  \\ \hline
	$3$  & $18$ & $18$ & $18$ & $18$ & $17$  \\ \hline
	$4$  & $20$ & $20$ & $20$ & $19$ & $19$  \\ \hline
	$5$  & $22$ & $22$ & $22$ & $21$ & $20$   \\
    \hline
  \end{tabular}
    \caption{Iteration numbers for Yeti-footprint using $\Pi_{\mathrm{cn}}$ and $M_{\mathrm{sD},1}$ (left), $M_{\mathrm{sD},2}$ (right)}
    \label{SognTakacs:tab:YPN}
  \end{center} 
\end{table}
\begin{figure}[tp]
		\begin{center}
				\resizebox{5cm}{5cm}{%
			\begin{tikzpicture}
				\begin{axis}[
					xlabel={Refinement level $\ell$ ($p=4$)},
					ylabel={Condition numbers},
					xmin=2, xmax=5,
					ymin=2, ymax=1000,
					xtick={2,3,4,5},
					ytick={2,5,10,20,50,100,200,500},
					legend pos=north west,
					ymajorgrids=true,
					grid style=dashed,
					legend columns=3,
					legend pos=north west,
					legend image post style={only marks},
					log ticks with fixed point,
					xmode=log,
					ymode=log,
					width=6cm
					]
					
					\addplot[color=blue,
					mark=triangle]
					coordinates {
						(2,114.87)(3,140.757)(4,171.456)(5,204.894)
					};
					\addlegendentry{{\tiny $M_{\mathrm{sD},1}$, $\Pi_{\mathrm{c}}$}};
					
					\addplot[color=red,
					mark=square]
					coordinates {
						(2,6.79309)(3,7.67078)(4,8.8891)(5,10.3149)
					};
					\addlegendentry{{\tiny $M_{\mathrm{sD},1}$, $\Pi_{\mathrm{ce}}$}};
					
					\addplot[color=green ,
					mark=o]
					coordinates {
						(2,14.8929)(3,18.082)(4,21.509)(5,25.2173)
					};
					\addlegendentry{{\tiny $M_{\mathrm{sD},1}$, $\Pi_{\mathrm{cn}}$}};

					\addplot[color=blue,
					mark=triangle*, densely dotted]
					coordinates {
						(2,41.207)(3,49.7688)(4,59.1471)(5,69.4346)
					};
					\addlegendentry{{\tiny $M_{\mathrm{sD},2}$, $\Pi_{\mathrm{c}}$}};
					
					\addplot[color=red,
					mark=square*, densely dotted]
					coordinates {
						(2,3.90832)(3,4.61709)(4,5.5316)(5,6.54666)
					};
					\addlegendentry{{\tiny $M_{\mathrm{sD},2}$, $\Pi_{\mathrm{ce}}$}};
					
					\addplot[color=green ,
					mark=*, densely dotted]
					coordinates {
						(2,8.53874)(3,10.3001)(4,11.6639)(5,14.2968)

					};
					\addlegendentry{{\tiny $M_{\mathrm{sD},2}$, $\Pi_{\mathrm{cn}}$}};
					
				\end{axis}
			\end{tikzpicture}
		}
		\quad
		\resizebox{5cm}{5cm}{%
			\begin{tikzpicture}
				\begin{axis}[
					xlabel={Polynomial degree $p$ ($\ell=5$)},
					ylabel={Condition numbers},
					xmin=2, xmax=6,
					ymin=5, ymax=1000,
					xtick={2,3,4,5,6},
					ytick={5,10,20,50,100,200,500},
					legend pos=north west,
					ymajorgrids=true,
					grid style=dashed,
					legend columns=3,
					legend pos=north west,
					legend image post style={only marks},
					log ticks with fixed point,
					xmode=log,
					ymode=log,
					width=6cm
					]
					
					\addplot[color=blue,
					mark=triangle]
					coordinates {
						(2,172.757)	(3,192.404)	(4,204.894)	(5,216.796)	(6,229.049)
					};
					\addlegendentry{{\tiny $M_{\mathrm{sD},1}$, $\Pi_{\mathrm{c}}$}};
					
					\addplot[color=red,
					mark=square]
					coordinates {
						(2,9.06913)	(3,9.73774)	(4,10.3149)	(5,10.8264)	(6,11.2882)
					};
					\addlegendentry{{\tiny $M_{\mathrm{sD},1}$, $\Pi_{\mathrm{ce}}$}};
					
					\addplot[color=green ,
					mark=o]
					coordinates {
						(2,21.7415)	(3,23.6535)	(4,25.2173)	(5,26.5506)	(6,27.7153)
					};
					\addlegendentry{{\tiny $M_{\mathrm{sD},1}$, $\Pi_{\mathrm{cn}}$}};

					\addplot[color=blue,
					mark=triangle*, densely dotted]
					coordinates {
						(2,60.0606)	(3,65.2032)	(4,69.4346)	(5,73.0702)	(6,76.267)
					};
					\addlegendentry{{\tiny $M_{\mathrm{sD},2}$, $\Pi_{\mathrm{c}}$}};
					
					\addplot[color=red,
					mark=square*, densely dotted]
					coordinates {
						(2,5.65636)	(3,5.95443)	(4,6.54666)	(5,6.7984)	(6,6.88155)
					};
					\addlegendentry{{\tiny $M_{\mathrm{sD},2}$, $\Pi_{\mathrm{ce}}$}};
					
					\addplot[color=green ,
					mark=*, densely dotted]
					coordinates {
						(2,12.7036)	(3,13.7621)	(4,14.2968)	(5,14.7683)	(6,15.5961)
					};
					\addlegendentry{{\tiny $M_{\mathrm{sD},2}$, $\Pi_{\mathrm{cn}}$}};
					
				\end{axis}
			\end{tikzpicture}
		}
	\end{center}
	\caption{Condition numbers for the quarter annulus\label{SognTakacs:figA}}
\end{figure}
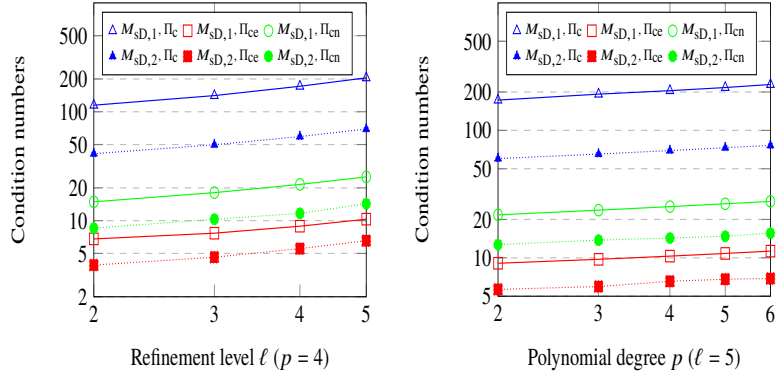

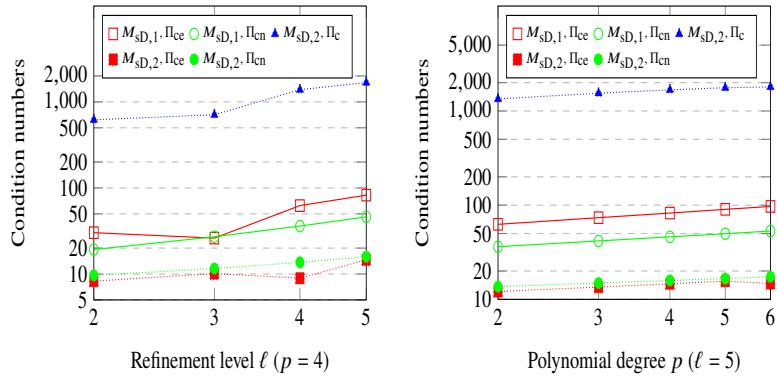
\begin{figure}[tp]
		\begin{center}
				\resizebox{5cm}{5cm}{%
			\begin{tikzpicture}
				\begin{axis}[
					xlabel={Refinement level $\ell$ ($p=4$)},
					ylabel={Condition numbers},
					xmin=2, xmax=5,
					ymin=5, ymax=13000,
					xtick={2,3,4,5},
					ytick={5,10,20,50,100,200,500,1000,2000},
					legend pos=north west,
					ymajorgrids=true,
					grid style=dashed,
					legend columns=3,
					legend pos=north west,
					legend image post style={only marks},
					log ticks with fixed point,
					xmode=log,
					ymode=log,
					width=6cm
					]

					\addplot[color=red,
					mark=square]
					coordinates {
						(2,30.2048)(3,26.1205)(4,62.3341)(5,82.4614)
					};
					\addlegendentry{{\tiny $M_{\mathrm{sD},1}$, $\Pi_{\mathrm{ce}}$}};
					
					\addplot[color=green ,
					mark=o]
					coordinates {
						(2,19.2215)(3,27.103)(4,36.0714)(5,46.1186)
					};
					\addlegendentry{{\tiny $M_{\mathrm{sD},1}$, $\Pi_{\mathrm{cn}}$}};

					\addplot[color=blue,
					mark=triangle*, densely dotted]
					coordinates {
						(2,619.296)(3,709.472)(4,1385.12)(5,1678.25)
					};
					\addlegendentry{{\tiny $M_{\mathrm{sD},2}$, $\Pi_{\mathrm{c}}$}};
					
					\addplot[color=red,
					mark=square*, densely dotted]
					coordinates {
						(2,8.26555)(3,10.1181)(4,8.93856)(5,14.644)
					};
					\addlegendentry{{\tiny $M_{\mathrm{sD},2}$, $\Pi_{\mathrm{ce}}$}};
					
					\addplot[color=green ,
					mark=*, densely dotted]
					coordinates {
						(2,9.65125)(3,11.5828)(4,13.6663)(5,15.8818)
					};
					\addlegendentry{{\tiny $M_{\mathrm{sD},2}$, $\Pi_{\mathrm{cn}}$}};
					
				\end{axis}
			\end{tikzpicture}
		}
		\quad
		\resizebox{5cm}{5cm}{%
			\begin{tikzpicture}
				\begin{axis}[
					xlabel={Polynomial degree $p$ ($\ell=5$)},
					ylabel={Condition numbers},
					xmin=2, xmax=6,
					ymin=10, ymax=13000,
					xtick={2,3,4,5,6},
					ytick={10,20,50,100,200,500,1000,2000,5000},
					legend pos=north west,
					ymajorgrids=true,
					grid style=dashed,
					legend columns=3,
					legend pos=north west,
					legend image post style={only marks},
					log ticks with fixed point,
					xmode=log,
					ymode=log,
					width=6cm
					]

					\addplot[color=red,
					mark=square]
					coordinates {
						(2,62.6185)	(3,73.6714)	(4,82.4614)	(5,90.3392)	(6,96.918)
					};
					\addlegendentry{{\tiny $M_{\mathrm{sD},1}$, $\Pi_{\mathrm{ce}}$}};
					
					\addplot[color=green ,
					mark=o]
					coordinates {
						(2,36.236)	(3,41.7154)	(4,46.1186)	(5,49.7836)	(6,53.1281)
					};
					\addlegendentry{{\tiny $M_{\mathrm{sD},1}$, $\Pi_{\mathrm{cn}}$}};

					\addplot[color=blue,
					mark=triangle*, densely dotted]
					coordinates {
						(2,1342.52)	(3,1543.77)	(4,1678.25)	(5,1764.14)	(6,1802.07)
					};
					\addlegendentry{{\tiny $M_{\mathrm{sD},2}$, $\Pi_{\mathrm{c}}$}};
					
					\addplot[color=red,
					mark=square*, densely dotted]
					coordinates {
						(2,12.1187)	(3,13.4986)	(4,14.644)	(5,15.6119)	(6,14.7552)
					};
					\addlegendentry{{\tiny $M_{\mathrm{sD},2}$, $\Pi_{\mathrm{ce}}$}};
					
					\addplot[color=green ,
					mark=*, densely dotted]
					coordinates {
						(2,13.6792)	(3,14.9083)	(4,15.8818)	(5,16.6108)	(6,17.3896)
					};
					\addlegendentry{{\tiny $M_{\mathrm{sD},2}$, $\Pi_{\mathrm{cn}}$}};
					
				\end{axis}
			\end{tikzpicture}
		}
	\end{center}
	\caption{Condition numbers for the Yeti-footprint\label{SognTakacs:figY}} 
\end{figure}

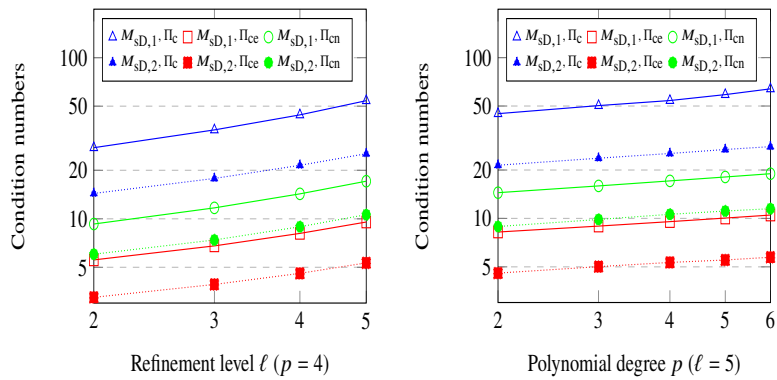
\begin{figure}[tp]
		\begin{center}
				\resizebox{5cm}{5cm}{%
			\begin{tikzpicture}
				\begin{axis}[
					xlabel={Refinement level $\ell$ ($p=4$)},
					ylabel={Condition numbers},
					xmin=2, xmax=5,
					ymin=3, ymax=200,
					xtick={2,3,4,5},
					ytick={5,10,20,50,100},
					legend pos=north west,
					ymajorgrids=true,
					grid style=dashed,
					legend columns=3,
					legend pos=north west,
					legend image post style={only marks},
					log ticks with fixed point,
					xmode=log,
					ymode=log,
					width=6cm
					]
					
					\addplot[color=blue,
					mark=triangle]
					coordinates {
						(2,27.6317)(3,35.6262)(4,44.1542)(5,54.0877)
					};
					\addlegendentry{{\tiny $M_{\mathrm{sD},1}$, $\Pi_{\mathrm{c}}$}};
					
					\addplot[color=red,
					mark=square]
					coordinates {
						(2,5.55832)(3,6.7911)(4,8.11693)(5,9.55116)
					};
					\addlegendentry{{\tiny $M_{\mathrm{sD},1}$, $\Pi_{\mathrm{ce}}$}};
					
					\addplot[color=green ,
					mark=o]
					coordinates {
						(2,9.28402)(3,11.6908)(4,14.2842)(5,17.1175)
					};
					\addlegendentry{{\tiny $M_{\mathrm{sD},1}$, $\Pi_{\mathrm{cn}}$}};

					\addplot[color=blue,
					mark=triangle*, densely dotted]
					coordinates {
						(2,14.3739)(3,17.8053)(4,21.4793)(5,25.3442)
					};
					\addlegendentry{{\tiny $M_{\mathrm{sD},2}$, $\Pi_{\mathrm{c}}$}};
					
					\addplot[color=red,
					mark=square*, densely dotted]
					coordinates {
						(2,3.25594)(3,3.91571)(4,4.58545)(5,5.32687)
					};
					\addlegendentry{{\tiny $M_{\mathrm{sD},2}$, $\Pi_{\mathrm{ce}}$}};
					
					\addplot[color=green ,
					mark=*, densely dotted]
					coordinates {
						(2,6.03108)(3,7.39175)(4,8.92729)(5,10.6089)
					};
					\addlegendentry{{\tiny $M_{\mathrm{sD},2}$, $\Pi_{\mathrm{cn}}$}};
					
				\end{axis}
			\end{tikzpicture}
		}
		\quad
		\resizebox{5cm}{5cm}{%
			\begin{tikzpicture}
				\begin{axis}[
					xlabel={Polynomial degree $p$ ($\ell=5$)},
					ylabel={Condition numbers},
					xmin=2, xmax=6,
					ymin=3, ymax=200,
					xtick={2,3,4,5,6},
					ytick={5,10,20,50,100},
					legend pos=north west,
					ymajorgrids=true,
					grid style=dashed,
					legend columns=3,
					legend pos=north west,
					legend image post style={only marks},
					log ticks with fixed point,
					xmode=log,
					ymode=log,
					width=6cm
					]
					
					\addplot[color=blue,
					mark=triangle]
					coordinates {
						(2,44.8186)	(3,50.3068)	(4,54.0877)	(5,58.882)	(6,63.9428)
					};
					\addlegendentry{{\tiny $M_{\mathrm{sD},1}$, $\Pi_{\mathrm{c}}$}};
					
					\addplot[color=red,
					mark=square]
					coordinates {
						(2,8.23096)	(3,8.96144)	(4,9.55116)	(5,10.0536)	(6,10.4935)
					};
					\addlegendentry{{\tiny $M_{\mathrm{sD},1}$, $\Pi_{\mathrm{ce}}$}};
					
					\addplot[color=green ,
					mark=o]
					coordinates {
						(2,14.4594)	(3,15.9295)	(4,17.1175)	(5,18.1048)	(6,18.9687)
					};
					\addlegendentry{{\tiny $M_{\mathrm{sD},1}$, $\Pi_{\mathrm{cn}}$}};

					\addplot[color=blue,
					mark=triangle*, densely dotted]
					coordinates {
						(2,21.4163)	(3,23.6869)	(4,25.3442)	(5,26.8281)	(6,27.9546)
					};
					\addlegendentry{{\tiny $M_{\mathrm{sD},2}$, $\Pi_{\mathrm{c}}$}};
					
					\addplot[color=red,
					mark=square*, densely dotted]
					coordinates {
						(2,4.57048)	(3,5.01424)	(4,5.32687)	(5,5.5146)	(6,5.73855)
					};
					\addlegendentry{{\tiny $M_{\mathrm{sD},2}$, $\Pi_{\mathrm{ce}}$}};
					
					\addplot[color=green ,
					mark=*, densely dotted]
					coordinates {
						(2,8.90624)	(3,9.85481)	(4,10.6089)	(5,11.116)	(6,11.4959)
					};
					\addlegendentry{{\tiny $M_{\mathrm{sD},2}$, $\Pi_{\mathrm{cn}}$}};
					
				\end{axis}
			\end{tikzpicture}
		}
	\end{center}
	\caption{Condition numbers for the unit square\label{SognTakacs:figU}} 
\end{figure}

In the Tables~\ref{SognTakacs:tab:APC}, \ref{SognTakacs:tab:APE} and \ref{SognTakacs:tab:APN}, we present the results for the quarter annulus domain. In these tables, we indicate the number of iterations of the preconditioned conjugate gradient solver required to reach the desired threshold. In the tables on the left side, we present the results obtained using the Stokes based scaled Dirichlet preconditioner $M_{\mathrm{sD},1}$. The tables on the right side show the results obtained using the Poisson based scaled Dirichlet preconditioner $M_{\mathrm{sD},2}$. The Tables~\ref{SognTakacs:tab:APC}, \ref{SognTakacs:tab:APE} and \ref{SognTakacs:tab:APN} show the results for the primal degrees of freedom $\Pi_{\mathrm{c}}$, $\Pi_{\mathrm{ce}}$ and $\Pi_{\mathrm{cn}}$, respectively. We observe convergence for all setups. We observe that only choosing the corner values as primal degrees of freedom (variant $\Pi_{\mathrm{c}}$) leads to the largest iteration numbers. The variant $\Pi_{\mathrm{ce}}$, which has to a slightly larger primal problem than $\Pi_{\mathrm{cn}}$, leads to the smallest iteration numbers.

We observe that the preconditioner $M_{\mathrm{sD},2}$ outperforms the preconditioner $M_{\mathrm{sD},1}$ in all cases. Note that the local problems required to be solved to realize $M_{\mathrm{sD},2}$ are smaller and they involve symmetric and positive definite matrices. So, $M_{\mathrm{sD},2}$ can also be realized more efficiently than $M_{\mathrm{sD},1}$.

We have estimated the condition numbers of the preconditioned system with the conjugate gradient solver. We present these results in Figure~\ref{SognTakacs:figA}. The left diagram shows the dependence on the refinement level $\ell$ for a fixed choice of the spline degree $p=4$. In the right diagram, one can see the dependence on the spline degree $p$ for a fixed choice of the refinement level $\ell=5$. In both cases, one can observe a mild dependence. The convergence of the solver for the variants $\Pi_{\mathrm{ce}}$ and $\Pi_{\mathrm{cn}}$ is close to what one would estimate from the condition number, that is, the convergence rate is only slightly smaller than
\[
		\frac{\sqrt{\kappa}-1}{\sqrt{\kappa}+1},
\]
where $\kappa$ is the condition number of the preconditioned system. If we consider the variant $\Pi_{\mathrm{c}}$, the difference is much more significant, that is, the true iteration numbers are much smaller than the estimated condition number would suggest. 

The next domain to consider is the Yeti-footprint.
Note that the inf-sup stability highly depends on the shape of the domain. Long and thin channels, which constitute the Yeti-footprint, are known to lead to small inf-sup constants.
In the Tables~\ref{SognTakacs:tab:YPC}, \ref{SognTakacs:tab:YPE} and \ref{SognTakacs:tab:YPN} and in Figure~\ref{SognTakacs:figY}, we present the results for the Yeti-footprint. The iteration counts and the condition numbers are only slightly larger than for the quarter annulus.  The largest difference is obtained for the choice $\Pi_{\mathrm{c}}$, particularly if the condition number is considered. If the preconditioner $M_{\mathrm{sD},1}$ is used, the condition number estimates gave a ``not a number'' results in many cases. Consequently, they were not added to the diagram.  Concerning the dependence on the grid size and the spline degree, the choice of the primal degrees of freedom and the choice of the preconditioner, we observe qualitatively the same results as for the quarter annulus.

Finally, we present results for the unit square in order to better see the effect of the geometry transformation. In the Tables~\ref{SognTakacs:tab:UPC}, \ref{SognTakacs:tab:UPE} and \ref{SognTakacs:tab:UPN} and in Figure~\ref{SognTakacs:figU}, one can see the results for this domain. In any case, the iteration counts and the condition numbers are smaller than those for the non-trivial domains. This difference is only mild for the choices $\Pi_{\mathrm{ce}}$ and $\Pi_{\mathrm{cn}}$, however for the choice $\Pi_{\mathrm{c}}$, this difference is quite large. When one only compares the results for the unit square, one would conclude that choice $\Pi_{\mathrm{c}}$ would be acceptable as well. The results for the other computational domains however show that $\Pi_{\mathrm{c}}$ is significantly inferior.

\section{Conclusions and final remarks}
\label{SognTakacs:sec:conclusions}
We solved the Stokes equations, discretized using multipatch IgA, by means of ITEI-DP solvers. Even though, the Stokes system is indefinite, the reduced problem is symmetric positive definite, which we can solve efficiently using a preconditioned conjugate gradient solver. Two scaled Dirichlet preconditioners were tested. The simpler one, which is based on the vector valued Poisson problem and thus easier to realize in practice, is numerically superior. Concerning the choice of the primal degrees of freedom, we observe that it is worthwhile to include edge averages of the velocity value (either for both components or only for the normal component). This observation cannot be made if only the unit square is considered. Thus, it is necessary to test ITEI-DP methods on non-trivial domains. 

Convergence analysis will be discussed in a forthcoming paper.

\begin{table}[tp]
  \begin{center} 
  \begin{tabular}{|c|c|c|c|c|c|}
    \hline
    $\;\ell\setminus p\;$ & {\quad$2$\quad} & {\quad$3$\quad} & {\quad$4$\quad} & {\quad$5$\quad} & {\quad$6$\quad} \\
    \hline
    \hline
  $2$  & $28$ & $29$ & $30$ & $28$ & $27$  \\ \hline
	$3$  & $32$ & $32$ & $33$ & $32$ & $31$  \\ \hline
	$4$  & $37$ & $37$ & $37$ & $36$ & $35$  \\ \hline
	$5$  & $40$ & $43$ & $42$ & $42$ & $39$   \\
    \hline
  \end{tabular}
  $\quad$
    \begin{tabular}{|c|c|c|c|c|c|}
    \hline
    $\;\ell\setminus p\;$ & {\quad$2$\quad} & {\quad$3$\quad} & {\quad$4$\quad} & {\quad$5$\quad} & {\quad$6$\quad} \\
    \hline
    \hline
  $2$  & $21$ & $21$ & $22$ & $21$ & $20$  \\ \hline
	$3$  & $23$ & $24$ & $24$ & $24$ & $23$  \\ \hline
	$4$  & $27$ & $26$ & $27$ & $26$ & $25$  \\ \hline
	$5$  & $27$ & $30$ & $30$ & $29$ & $28$  \\
    \hline
  \end{tabular}
    \caption{Iteration numbers for unit square using $\Pi_{\mathrm{c}}$ and $M_{\mathrm{sD},1}$ (left), $M_{\mathrm{sD},2}$ (right)}
    \label{SognTakacs:tab:UPC}
  \end{center} 
\end{table}
\begin{table}[tp]
  \begin{center} 
  \begin{tabular}{|c|c|c|c|c|c|}
    \hline
    $\;\ell\setminus p\;$ & {\quad$2$\quad} & {\quad$3$\quad} & {\quad$4$\quad} & {\quad$5$\quad} & {\quad$6$\quad} \\
    \hline
    \hline
  $2$  & $14$ & $14$ & $14$ & $14$ & $14$  \\ \hline
	$3$  & $15$ & $16$ & $16$ & $16$ & $15$  \\ \hline
	$4$  & $17$ & $17$ & $17$ & $17$ & $17$  \\ \hline
	$5$  & $18$ & $18$ & $18$ & $18$ & $18$  \\
    \hline
  \end{tabular}
  $\quad$
    \begin{tabular}{|c|c|c|c|c|c|}
    \hline
    $\;\ell\setminus p\;$ & {\quad$2$\quad} & {\quad$3$\quad} & {\quad$4$\quad} & {\quad$5$\quad} & {\quad$6$\quad} \\
    \hline
    \hline
  $2$  & $18$ & $18$ & $18$ & $19$ & $18$  \\ \hline
	$3$  & $20$ & $21$ & $21$ & $21$ & $20$  \\ \hline
	$4$  & $22$ & $23$ & $22$ & $22$ & $22$  \\ \hline
	$5$  & $24$ & $25$ & $25$ & $24$ & $23$  \\
    \hline
  \end{tabular}
    \caption{Iteration numbers for unit square using $\Pi_{\mathrm{ce}}$ and $M_{\mathrm{sD},1}$ (left), $M_{\mathrm{sD},2}$ (right)}
    \label{SognTakacs:tab:UPE}
  \end{center} 
\end{table}
\begin{table}[tp]
  \begin{center} 
  \begin{tabular}{|c|c|c|c|c|c|}
    \hline
    $\;\ell\setminus p\;$ & {\quad$2$\quad} & {\quad$3$\quad} & {\quad$4$\quad} & {\quad$5$\quad} & {\quad$6$\quad} \\
    \hline
    \hline
  $2$  & $18$ & $18$ & $18$ & $19$ & $18$  \\ \hline
	$3$  & $20$ & $21$ & $21$ & $21$ & $20$  \\ \hline
	$4$  & $22$ & $23$ & $22$ & $22$ & $22$  \\ \hline
	$5$  & $24$ & $25$ & $25$ & $24$ & $23$   \\
    \hline
  \end{tabular}
  $\quad$
    \begin{tabular}{|c|c|c|c|c|c|}
    \hline
    $\;\ell\setminus p\;$ & {\quad$2$\quad} & {\quad$3$\quad} & {\quad$4$\quad} & {\quad$5$\quad} & {\quad$6$\quad} \\
    \hline
    \hline
	$2$  & $14$ & $15$ & $15$ & $15$ & $14$  \\ \hline
	$3$  & $16$ & $16$ & $16$ & $16$ & $16$  \\ \hline
	$4$  & $18$ & $18$ & $18$ & $18$ & $17$  \\ \hline
	$5$  & $19$ & $20$ & $20$ & $20$ & $19$  \\
    \hline
  \end{tabular}
    \caption{Iteration numbers for unit square using $\Pi_{\mathrm{cn}}$ and $M_{\mathrm{sD},1}$ (left), $M_{\mathrm{sD},2}$ (right)}
    \label{SognTakacs:tab:UPN}
  \end{center} 
\end{table}

\begin{acknowledgement}
This work was supported by the Austrian Science Fund (FWF): P31048. This support is gratefully acknowledged.
\end{acknowledgement}

\bibliographystyle{spmpsci}
\bibliography{bibliography}
\end{document}